\expandafter\ifx\csname mthreemacsloaded\endcsname\relax\else \fi

\magnification1100
\input amstex


 \catcode`\@=11
 \let\wlog@ld\wlog
 \def\wlog#1{\relax}

 \newif\ifIN@
 \def\m@rker{\m@@rker}
 \def\IN@{\expandafter\INN@\expandafter}
 \long\def\INN@0#1@#2@{\long\def\NI@##1#1##2##3\ENDNI@
    {\ifx\m@rker##2\IN@false\else\IN@true\fi}%
     \expandafter\NI@#2@@#1\m@rker\ENDNI@}
  \newtoks\Initialtoks@  \newtoks\Terminaltoks@
  \def\SPLIT@{\expandafter\SPLITT@\expandafter}
  \def\SPLITT@0#1@#2@{\def\TTILPS@##1#1##2@{%
     \Initialtoks@{##1}\Terminaltoks@{##2}}\expandafter\TTILPS@#2@}
  \newtoks\Trimtoks@

 \def\ForeTrim@{\expandafter\ForeTrim@@\expandafter}
 \def\ForePrim@0 #1@{\Trimtoks@{#1}}
 \def\ForeTrim@@0#1@{\IN@0\m@rker. @\m@rker.#1@%
     \ifIN@\ForePrim@0#1@%
     \else\Trimtoks@\expandafter{#1}\fi}
 
  \def\Trim@0#1@{%
      \ForeTrim@0#1@%
      \IN@0 @\the\Trimtoks@ @%
        \ifIN@
             \SPLIT@0 @\the\Trimtoks@ @\Trimtoks@\Initialtoks@
             \IN@0\the\Terminaltoks@ @ @%
                 \ifIN@
                 \else \Trimtoks@ {FigNameWithSpace}%
                 \fi
        \fi
      }

  \font\titlebold=cmbx12 scaled 1200
  \font\twelvebold=cmbx12
  \font\tenbold=cmbx10
  \font\ninebold=cmbx9
  \font\sevenbold=cmbx7
  \font\fivebold=cmbx5

  \input amssym.def \input amssym
     \font\titlemsa=msam10 at 14.4pt
     \font\titlemsb=msbm10 at 14.4pt
     \font\titleeufm=eufm10 at 14.4pt
     \font\twelvemsa=msam10 scaled 1200
     \font\twelvemsb=msbm10 scaled 1200
     \font\twelveeufm=eufm10 scaled 1200
     \font\ninemsa=msam9
     \font\ninemsb=msbm9
     \font\nineeufm=eufm9

   \ifx\cyrfam\undefined
   \else
     \immediate\write16{}%
     \message{ !!! cyr fonts already defined. !!! }
     \message{ --- edit out superfluous font defs? }
   \fi
   \newfam\cyrfam
       \font\titlecyr=wncyr10 scaled 1440 
       \font\twelvecyr=wncyr10 scaled 1200
       \font\tencyr=wncyr10
       \font\ninecyr=wncyr9
       \font\sevencyr=wncyr7
       \font\sixcyr=wncyr6

   \newfam\eusmfam
       \font\titleeusm=eusm10 scaled 1440
       \font\twelveeusm=eusm10 scaled 1200
       \font\teneusm=eusm10
       \font\nineeusm=eusm9
       \font\seveneusm=eusm7
       
       \font\fiveeusm=eusm5

\let\Cal\cal

    \font\ninemrm=cmr9 
    \font\ninei=cmmi9
    \font\ninesy=cmsy9 
    \skewchar\ninei='177
    \skewchar\ninesy='60

  \font\twelvemrm=cmr10 at 12pt 
  \font\twelvei=cmmi10 at 12pt
  \font\twelvesy=cmsy10 at 12pt

  \font\titlemrm=cmr10 at 14.4pt 
  \font\titlei=cmmi10 at 14.4pt
  \font\titlesy=cmsy10 at 14.4pt


  \def\Smallfonts{\ninepoint}

  \def\Hfont{\titlepoint\bf}
  \def\Authorfont{\twelvepoint\it}
  \def\HHfont{\twelvepoint\bf}
  \def\HHHfont{\bf}
  \def\Bibfont{\tenbf}
  \def\Coordfont{\nineit }

  \def \thfont {\bf }
  \def \pffont {\it\itSpacing }
  \def \rkfont {\bf }
  \def \dffont {\bf }
  \def \egfont {\bf }

 \def\ninepoint{%
  \def\rm{\fam0\ninerm}%
    \textfont0=\ninemrm  \scriptfont0=\sevenrm  \scriptscriptfont0=\fiverm
    \textfont1=\ninei    \scriptfont1=\seveni   \scriptscriptfont1=\fivei
  \def\mit{\fam1\ninei}%
  \def\oldstyle{\fam1\ninei}%
    \textfont2=\ninesy   \scriptfont2=\sevensy  \scriptscriptfont2=\fivesy
    \textfont3=\tenex    \scriptfont3=\tenex    \scriptscriptfont3=\tenex
  \def\it{\fam\itfam\nineit}%
    \textfont\itfam=\nineit
  \def\bf{\ifmmode\fam\bffam\else\ninebf\fi}%
    \textfont\bffam=\ninebold 
    \scriptfont\bffam=\sevenbold 
    \scriptscriptfont\bffam=\fivebold%
  \def\msa{\fam\msafam\ninemsa}%
    \textfont\msafam=\ninemsa 
    \scriptfont\msafam=\sevenmsa
    \scriptscriptfont\msafam=\fivemsa%
  \def\msb{\fam\msbfam\ninemsb}%
    \textfont\msbfam=\ninemsb%
    \scriptfont\msbfam=\sevenmsb%
    \scriptscriptfont\msbfam=\fivemsb%
  \def\eufm{\fam\eufmfam\nineeufm}%
    \textfont\eufmfam=\nineeufm
    \scriptfont\eufmfam=\seveneufm
    \scriptscriptfont\eufmfam=\fiveeufm
   \def\eusm{\fam\eusmfam\nineeusm}%
     \textfont\eusmfam=\nineeusm
     \scriptfont\eusmfam=\seveneusm
     \scriptscriptfont\eusmfam=\fiveeusm
   \def\cyr{\fam\cyrfam\ninecyr}%
     \textfont\cyrfam=\ninecyr
     \scriptfont\cyrfam=\sevencyr
     \scriptscriptfont\cyrfam=\sixcyr
  \setbox\strutbox=\hbox{\vrule
      height7pt depth3pt width0pt}%
   \baselineskip=10.8pt\rm}

 \let\eightpoint\ninepoint 

 \def\tenpoint{%
  \def\rm{\fam0\tenrm}%
    \textfont0=\tenmrm \scriptfont0=\sevenrm \scriptscriptfont0=\fiverm%
  \def\mit{\fam1\teni}%
  \def\oldstyle{\fam1\teni}%
    \textfont1=\teni   \scriptfont1=\seveni  \scriptscriptfont1=\fivei%
    \textfont2=\tensy  \scriptfont2=\sevensy \scriptscriptfont2=\fivesy%
    \textfont3=\tenex  \scriptfont3=\tenex   \scriptscriptfont3=\tenex%
  \def\it{\fam\itfam\tenit}%
    \textfont\itfam=\tenit%
  \def\bf{\ifmmode\fam\bffam\else\tenbf\fi}%
    \textfont\bffam=\tenbold
    \scriptfont\bffam=\sevenbold%
    \scriptscriptfont\bffam=\fivebold%
  \def\msa{\fam\msafam\tenmsa}%
    \textfont\msafam=\tenmsa%
    \scriptfont\msafam=\sevenmsa%
    \scriptscriptfont\msafam=\fivemsa%
  \def\msb{\fam\msbfam\tenmsb}%
    \textfont\msbfam=\tenmsb%
    \scriptfont\msbfam=\sevenmsb%
    \scriptscriptfont\msbfam=\fivemsb%
  \def\eufm{\fam\eufmfam\teneufm}%
   \textfont\eufmfam=\teneufm
   \scriptfont\eufmfam=\seveneufm
   \scriptscriptfont\eufmfam=\fiveeufm
   \def\eusm{\fam\eusmfam\teneusm}%
    \textfont\eusmfam=\teneusm
    \scriptfont\eusmfam=\seveneusm
    \scriptscriptfont\eusmfam=\fiveeusm
   \def\cyr{\fam\cyrfam\tencyr}%
    \textfont\cyrfam=\tencyr
    \scriptfont\cyrfam=\sevencyr
    \scriptscriptfont\cyrfam=\sixcyr
  \setbox\strutbox=\hbox{\vrule %
      height8.5pt depth3.5ptwidth0pt}%
  \baselineskip=\StdBaselineskip\rm}

 \def\twelvepoint{%
  \def\rm{\fam0\twelverm}%
    \textfont0=\twelvemrm \scriptfont0=\tenmrm \scriptscriptfont0=\sevenrm
    \textfont1=\twelvei   \scriptfont1=\teni   \scriptscriptfont1=\seveni
  \def\mit{\fam1\twelvei}%
  \def\oldstyle{\fam1\twelvei}%
    \textfont2=\twelvesy  \scriptfont2=\tensy  \scriptscriptfont2=\sevensy
    \textfont3=\tenex  \scriptfont3=\tenex  \scriptscriptfont3=\tenex
  \def\it{\fam\itfam\twelveit}%
    \textfont\itfam=\twelveit
  \def\bf{\ifmmode\fam\bffam\else\twelvebf\fi}%
    \textfont\bffam=\twelvebold
    \scriptfont\bffam=\tenbold%
    \scriptscriptfont\bffam=\sevenbold%
  \def\msa{\fam\msafam\twelvemsa}%
    \textfont\msafam=\twelvemsa%
    \scriptfont\msafam=\tenmsa%
    \scriptscriptfont\msafam=\sevenmsa%
  \def\msb{\fam\msbfam\twelvemsb}%
    \textfont\msbfam=\twelvemsb%
    \scriptfont\msbfam=\tenmsb%
    \scriptscriptfont\msbfam=\sevenmsb%
  \def\eufm{\fam\eufmfam\twelveeufm}%
   \textfont\eufmfam=\twelveeufm
   \scriptfont\eufmfam=\teneufm
   \scriptscriptfont\eufmfam=\seveneufm
   \def\eusm{\fam\eusmfam\twelveeusm}%
    \textfont\eusmfam=\twelveeusm
    \scriptfont\eusmfam=\teneusm
    \scriptscriptfont\eusmfam=\seveneusm
   \def\cyr{\fam\cyrfam\tencyr}%
    \textfont\cyrfam=\twelvecyr
    \scriptfont\cyrfam=\tencyr
    \scriptscriptfont\cyrfam=\sevencyr
  \setbox\strutbox=\hbox{\vrule
      height10.2pt depth4.55pt width0pt}%
  \baselineskip=14pt\rm}

 \def\titlepoint{%
    \textfont0=\titlemrm \scriptfont0=\twelvemrm \scriptscriptfont0=\tenmrm
    \textfont1=\titlei   \scriptfont1=\twelvei   \scriptscriptfont1=\teni
  \def\mit{\fam1\titlei}%
  \def\oldstyle{\fam1\titlei}%
    \textfont2=\titlesy  \scriptfont2=\twelvesy  \scriptscriptfont2=\tensy
    \textfont3=\tenex
    \scriptfont3=\tenex
    \scriptscriptfont3=\tenex
  \def\it{\fam\itfam\titleit}%
    \textfont\itfam=\titleit
  \def\bf{\ifmmode\fam\bffam\else\titlebf\fi}%
    \textfont\bffam=\titlebold
    \scriptfont\bffam=\twelvebold%
    \scriptscriptfont\bffam=\tenbold%
  \def\msa{\fam\msafam\titlemsa}%
    \textfont\msafam=\titlemsa%
    \scriptfont\msafam=\twelvemsa%
    \scriptscriptfont\msafam=\tenmsa%
  \def\msb{\fam\msbfam\titlemsb}%
    \textfont\msbfam=\titlemsb%
    \scriptfont\msbfam=\twelvemsb%
    \scriptscriptfont\msbfam=\tenmsb%
  \def\eufm{\fam\eufmfam\titleeufm}%
    \textfont\eufmfam=\titleeufm
    \scriptfont\eufmfam=\twelveeufm
    \scriptscriptfont\eufmfam=\teneufm
   \def\eusm{\fam\eusmfam\titleeusm}%
     \textfont\eusmfam=\titleeusm
     \scriptfont\eusmfam=\twelveeusm
     \scriptscriptfont\eusmfam=\teneusm
   \def\cyr{\fam\cyrfam\tencyr}%
    \textfont\cyrfam=\titlecyr
    \scriptfont\cyrfam=\twelvecyr
    \scriptscriptfont\cyrfam=\tencyr
  \setbox\strutbox=\hbox{\vrule
      height12.3pt depth5.54pt width0pt}%
  \baselineskip=16pt\rm}

\newbox\AuthorBox\newbox\TitleBox
\newbox\TFLinebox
\newbox\FLinebox
\newbox\HLinebox
\def\SetTFLinebox#1{\setbox\TFLinebox=\hbox{#1}}
\def\SetFLinebox#1{\setbox\FLinebox=\hbox{#1}}
\def\SetHLinebox#1{\setbox\HLinebox=\hbox{#1}}

 \def\SetAuthorHead#1{%
     \setbox\AuthorBox=\hbox{\ninepoint \it 
           \ignorespaces\frenchspacing#1\unskip}}
 \def\SetTitleHead#1{%
     \setbox\TitleBox=\hbox{\ninepoint \it
           \ignorespaces\frenchspacing#1\unskip}}

  \def\itSpacing{\relax}
  \def\itSpacingOff{\relax}


 \def\Hrule{\hrule width0pt height0pt}

  \newskip\ProcSkip \ProcSkip 8pt plus2pt minus2pt

 \newskip\LastSkip
 \def\SaveLastSkip{\LastSkip\lastskip}
 \def\RestoreLastSkip{\vskip-\LastSkip\vskip\LastSkip}

 \def\NoindentAfter{\everypar={\setbox0=\lastbox\everypar={}}}

 \long\def\H#1\par#2\par{\notenumber=0 \titlepagetrue%
    {
    \baselineskip=20pt
    \parindent=0pt\parskip=0pt\frenchspacing
    \leftskip=0pt plus .2\hsize minus .3\hsize
    \rightskip=0pt plus .2\hsize minus .3\hsize
 \def\\{\unskip\break}%
    \pretolerance=10000 \Hfont #1\unskip\break
     \vskip7pt\Hrule
\hfill \Authorfont #2\hfill\hfill\unskip}
    \vskip48pt plus 4pt minus 4pt
    \par\NoindentAfter\rm}

 \long\def\Hi#1\par#2\par{\notenumber=0 \titlepagetrue%
    {  \baselineskip=0pt  \parindent=0pt\parskip=0pt\frenchspacing
    \leftskip=0pt plus .2\hsize minus .3\hsize
    \rightskip=0pt plus .2\hsize minus .3\hsize
}
    \rm}


 \newdimen\PageRemainder
  \def\SetPageRemainder{
     \PageRemainder=\pagegoal
     \ifdim\PageRemainder=\maxdimen\PageRemainder=\vsize
     \else\advance\PageRemainder by -1\pagetotal\fi}

  \def\Rpt@{}\def\Rpt@@{}

  \long\def\HH#1\par{\par
  \SaveLastSkip\removelastskip\goodbreak
  \ifdim\LastSkip<30pt 
     \LastSkip 30pt
plus 3pt minus 2pt\fi
  \SetPageRemainder\advance\PageRemainder-\LastSkip
  \ifdim\PageRemainder<150pt
       \edef\Rpt@{remain = \the\PageRemainder\noexpand\\
                pagetotal=\the\pagetotal\noexpand\\
                           pagegoal=\the\pagegoal}%
          \fi
   \ifdim\PageRemainder<65pt 
       \ifdim\PageRemainder > 0pt
          \edef\Rpt@@{\noexpand\\
                      Had HH PageRemainder$<$\relax 65pt\noexpand\\
                      Hence forced break!}%
     \vskip 0pt plus .2\PageRemainder\eject 
    \fi\fi
    \vskip\LastSkip\Hrule 
    \pretolerance=10000\rightskip=0pt plus 3em
    \hangafter1 \hangindent=2.2em%
    \noindent
    \HHfont \unskip \Ednote{\Rpt@\Rpt@@}%
            \def\Rpt@{}\def\Rpt@@{}%
            \ignorespaces
            #1\par\rightskip=0pt\pretolerance=\StdPretolerance%
    \NoindentAfter
\tenpoint\rm%
     \medskip \vskip\ProcSkip}

  \long\def\HHH#1\par{\par%
  \SaveLastSkip\removelastskip\goodbreak
  \ifdim\LastSkip<\ProcSkip%
     \LastSkip\ProcSkip\fi
  \SetPageRemainder\advance\PageRemainder-\LastSkip
  \ifdim\PageRemainder<150pt
       \edef\Rpt@{remain = \the\PageRemainder\noexpand\\
                pagetotal=\the\pagetotal\noexpand\\
                           pagegoal=\the\pagegoal}%
       \fi
   \ifdim\PageRemainder<48pt  
        \ifdim\PageRemainder > 0pt
             \edef\Rpt@@{\noexpand\\
                      Had HHH PageRemainder$<$\relax48pt\noexpand\\
                      Hence forced break!}%
       \vskip 0pt plus .2\PageRemainder\eject 
      \fi\fi
   \vskip\LastSkip\par\noindent
   \HHHfont \unskip\Ednote{\Rpt@\Rpt@@}%
  \def\Rpt@{}\def\Rpt@@{}%
  \ignorespaces
   #1\unskip.\quad\rm\ignorespaces
   \ignorepars}

  \long\def\ignorepars#1\par{\def\Test{#1}%
     \ifx\Test\Empty\def\This{\ignorepars}%
        \else\def\This{\Test\par}\fi
           \This}
  \def\Empty{}

 \def\Abstract#1\par{\bgroup\Smallfonts\narrower\HHH #1\par}
 \def\endAbstract{\par\egroup}


 \def\ProcBreak{\par%
    \ifdim\lastskip<8pt%
    \removelastskip%
    \penalty-200\vskip\ProcSkip\fi}

 \def\th#1\par{\ProcBreak \noindent
   {\thfont\ignorespaces
    #1\unskip.}\it\itSpacing\kern.4em\ignorepars}

 \def\endth{\ProcBreak\rm\itSpacingOff }


 \def\pf#1\par{\ProcBreak %
    \noindent\pffont#1\unskip.\rm\itSpacingOff{\kern .7em}\ignorepars}


  \def\qedbox{\hbox{\vbox{
    \hrule width0.2cm height0.2pt
    \hbox to 0.2cm{\vrule height 0.2cm width 0.2pt
             \hfil\vrule height0.2cm width 0.2pt}
    \hrule width0.2cm height 0.2pt}\kern1pt}}

  \def\qed{\ifmmode\qedbox
    \else\unskip\ \hglue0mm\hfill\qedbox\ProcBreak\fi}

  \def \rk #1\par{\ProcBreak
     \noindent{\rkfont\ignorespaces #1\unskip.}%
     \rm\kern.6em\ignorepars}

  \def \df #1\par{\ProcBreak
     \noindent{\dffont\unskip\ignorespaces #1\unskip.}%
     \rm\kern.6em\ignorepars}

  \def \eg #1\par{\ProcBreak
     \noindent\egfont\unskip\ignorespaces #1\unskip.
     \rm\kern.6em\ignorepars}

  \newdimen\Overhang

   \def\MaxTag@#1#2#3#4#5{\setbox0=\hbox{#4\ignorespaces#2\unskip}%
     \dimen0=\wd0\advance\dimen0 by#3
     \ifdim\dimen0<#5\relax\dimen0=#5\fi
     \expandafter\edef\csname #1Hang\endcsname{\the\dimen0}}

 \def\MaxItemTag#1{\MaxTag@{Item}{#1}{.4em}{\ItemStyle}{\parindent}}%
 \def\MaxItemItemTag#1{%
        \MaxTag@{ItemItem}{#1}{.4em}{\ItemItemStyle}{\parindent}}
 \def\MaxNrTag#1{\MaxTag@{Nr}{#1}{.5em}{\NrStyle}{\parindent}}
 \def\MaxReferenceTag#1{%
        \MaxTag@{Reference}{[#1]}{.6em}{\ninerm}{\parindent}}
 \def\MaxFootTag#1{\MaxTag@{Foot}{#1}{.4em}{\ninerm}{\z@}}

  \def\SetOverhang@{\Overhang=.8\dimen0%
     \advance\Overhang by \wd0\relax
     \ifdim\Overhang>\hangindent\relax
       \advance\Overhang by .25\dimen0%
       \Ednote{Tag is pushing text.}\osumess{Tag is pushing text.}%
     \else\Overhang=\hangindent
     \fi}

   \def\Item#1{\par\noindent
      \hangafter1\hangindent=\ItemHang
      \setbox0=\hbox{\ItemStyle\ignorespaces#1\unskip}%
      \dimen0=.4em\SetOverhang@
      \rlap{\box0}\kern\Overhang\ignorespaces}

   \def\ItemItem#1{\par\noindent
      \hangafter1\hangindent=\ItemItemHang
      \setbox0=\hbox{\ItemItemStyle\ignorespaces#1\unskip}%
      \dimen0=.4em\SetOverhang@
      \advance\hangindent by \ItemHang
      \kern\ItemHang\rlap{\box0}%
      \kern\Overhang\ignorespaces}

  \def\Nr#1{\par\noindent\hangindent=\NrHang 
    \setbox0=\hbox{\NrStyle\ignorespaces#1\unskip}%
    \dimen0=.5em\SetOverhang@
    \rlap{\box0}\kern\Overhang
    \hangindent=\z@\ignorespaces}

   \newskip\Rosterskip\Rosterskip 1pt plus1pt 
   \def\Roster{\par\ifdim\lastskip<\Rosterskip\removelastskip\vskip\Rosterskip\fi
    \bgroup}
   \def\endRoster{\par\global\edef\LastSkip@{\the\lastskip}\removelastskip
       \egroup\penalty-50\LastSkip\LastSkip@\relax
       \ifdim\LastSkip<\Rosterskip\LastSkip\Rosterskip\fi
       \vskip\LastSkip}




 \def\cite#1{
    \def\nextiii@##1,##2\end@{{\frenchspacing\rm 
      \lBr\ignorespaces##1\unskip{\rm,~\ignorespaces##2}\rBr}}%
    \IN@0,@#1@%
    \ifIN@\def\next{\nextiii@#1\end@}\else
    \def\next{{\rm\lBr#1\rBr}}\fi\next}


   \def \Bib#1\par{%
       \par\removelastskip\SetPageRemainder
       \ifdim\PageRemainder < 97pt
        \ifdim\PageRemainder > 0pt
        \vfill\eject
       \fi\fi
    \ProcBreak \par\begingroup\parskip=0 pt%
    \goodbreak \vskip 15 pt plus 10 pt
    \noindent\null\hfill\Bibfont
      \ignorespaces #1\unskip\hfill\null\par 
    \frenchspacing \Smallfonts\rm
    \parskip=2.5 pt plus 1 pt minus.5pt%
    \nobreak\vskip 12pt plus 2pt minus2pt\nobreak
    \leftskip=0 pt \baselineskip=10.5pt}

 \def\ReferenceTagSlide{0em}
  \def\ReferenceTagGap{.5em}

  \def \rf#1{\par\noindent
     \hangafter1\hangindent=\ReferenceHang      
     \setbox0=\hbox{\ninerm[\ignorespaces#1\unskip]}%
     \dimen0=\ReferenceTagGap\SetOverhang@
     \rlap{\kern\ReferenceTagSlide\box0}%
     \kern\Overhang\ignorespaces}

  \def\ref#1\par#2\par#3\par#4\par{%
     \rf{#1}#2\unskip,\ #3\unskip,\
     #4\unskip.}

  \def\endBib{\par\endgroup\vskip 12pt minus 6pt }


  \long\def\Coordinates#1\endCoordinates{
 {\par\vskip4pt\def\\{\unskip, }\Coordfont\baselineskip10.5pt\noindent#1}}

 \def\pagecontents{
  \gdef\Pagetot@l{\pagetotal}
  \ifvoid\TRMargIns\else
    \rlap{\kern\hsize\kern10pt\vbox to 0pt{%
         \box\TRMargIns\vss}}\fi
  \ifvoid\topins\else\unvbox\topins\fi
   \dimen@=\dp\@cclv \unvbox\@cclv 
   \ifvoid\footins\else 
     \vskip\skip\footins
     \footnoterule
     \unvbox\footins\fi
   \ifr@ggedbottom \kern-\dimen@ \vfil \fi}


 \newcount\Ht 

 \def \Acc{\expandafter } 

 \def\swthat{\raise -1.1 ex\hbox{\sam$\widehat{}$}}
 \def\swttilde{\raise -1.2 ex\hbox{\sam$\widetilde{}$}}
 \def \overdot{{\raise .2 ex \hbox to 0pt {\hss\bf\smash{.}\hss}}}
 \def \overcircle{{\raise .1 ex \hbox to 0pt
    {\sam$\eightpoint\scriptstyle\hss\circ\hss$}}}

 \def \Mathaccent#1#2{{\sam 
  \setbox4=\hbox{$\vphantom{#2}$}
  \Ht=\ht4 
  \setbox5=\hbox{${#1}$}
  \setbox6=\hbox{${#2}$}
  \setbox7=\hbox to .5\wd6{}
  \copy7\kern .1\Ht \raise\Ht sp\hbox{\copy5}\kern-.1\Ht
  \copy7\llap{\box6}
  }}

  \def\SwtCheck #1{
        \ifmmode \check{#1}%
                \else \v {#1}%
                \fi}

 \def\barpartial {%
   \kern .17 em
    \overline {\kern -.17 em\partial\kern-.03 em}%
    \kern .03 em}

 
  \def\Overline#1{\setbox1=\hbox{\sam ${#1}$}%
      \ifdim \wd1 > 6pt
    \kern .11 em
    \overline {\kern -.11 em#1\kern-.14 em}
    \kern .14 em
  \else
    \kern .03 em
    \overline {\kern -.03 em#1\kern-.04 em}
    \kern .04 em
  \fi}

 \def\SOverline#1{\setbox1=\hbox{\sam ${#1}$}%
      \ifdim \wd1 > 7pt
    \kern .22 em
    \overline {\kern -.22 em#1\kern-.09 em}%
    \kern .09 em
  \else
    \kern .10 em
    \overline {\kern -.10 em#1\kern-.04 em}%
    \kern .04 em
  \fi}


 \def\Underline#1{\setbox1=\hbox{\sam ${#1}$}%
      \ifdim \wd1 > 6pt
    \kern .11 em
    \underline {\kern -.11 em#1\kern-.14 em}
    \kern .14 em
  \else
    \kern .03 em
    \underline {\kern -.03 em#1\kern-.04 em}
    \kern .04 em
  \fi}

 \def\SUnderline#1{\setbox1=\hbox{\sam ${#1}$}%
      \ifdim \wd1 > 7pt
    \kern .04 em
    \underline {\kern -.04 em#1\kern-.2 em}%
    \kern .2 em
  \else
    \kern .0 em
    \underline {\kern -.0 em#1\kern-.15 em}%
    \kern .15 em
  \fi}


 \def \Blackbox
   {\leavevmode\hskip .3pt \vbox
   {\hrule height 5pt\hbox{\hskip 4.5pt}}\hskip .5pt}

 \def \XX{\Blackbox\kern.5pt\Blackbox} 

  \def\.{.\kern1pt}

    \def\Hyphen{\edef\this{\the\hyphenchar\font}%
          \hyphenchar\font=-1\char\this\hyphenchar\font=\this}

 \ifx\undefined\text
  \def\text#1{\hbox{\rm #1}}\fi 



   \everymath{}  

  \def\PassMath@@{\aftergroup\AfterMath@} 

 \let\PassMath@\PassMath@@

 \def\AfterMath@{\futurelet\next\AfterMathMole@}

 \def\AfterMathMole@{
      \ifcat\next\space
          \def\this{}
      \else
      \ifcat\next\egroup %
        \def\this{\osumess{Handset mathsurround?? ---(see dollar brace)}}%
      \else
      \def\this{\AAfterMath@}
      \fi\fi
      \this}

 \def\hyphen@{-}
 \def\paren@{)}
 \def\apostr@{'}

 \def\MSC#1{\kern-.8\mathsurround#1\kern.8\mathsurround}

 \def\AAfterMath@#1{\def\Next{#1}
    \IN@0\Next @,.;:!?\relax @%
    \ifIN@\def\this{\MSC{\Next}}%
    \else
    \ifx\Next\hyphen@\def\this{\futurelet\next\AfterHyphen@}%
    \else
    \ifx\Next\paren@\def\this{#1}%
    \else 
    \ifx\Next\apostr@\def\this{#1}%
    \else \def\this{\osumess{Handset mathsurround??}%
                 #1}\fi\fi\fi\fi
    \this}

 \def\AfterHyphen@#1{\def\Next{#1}%
   \ifx\Next\hyphen@\def\this{--}\else
   \ifcat\next\space%
   \def\this{\kern-\mathsurround\kern.05em- \Next}\else
   \def\this{\kern-\mathsurround\kern.05em\Hyphen\Next}\fi\fi\this}

 \def\sam{\mathsurround=\z@\let\PassMath@\relax}  %
 \def\mas{\mathsurround=\StdMathsurround\let\PassMath@\PassMath@@}
 
 \def\Mas{\mathsurround=\StdMathsurround
                \everymath{\PassMath@}\let\PassMath@\PassMath@@}

 \def\m@th{\mathsurround=\z@\everymath{}}

 \def\m@@th{\mathsurround=\z@\everymath={}\let\m@th\relax}

\def\underbar#1{$\setbox\z@\hbox{#1}\dp\z@\z@
      \m@th \underline{\box\z@}$\relax}

\def\mathhexbox#1#2#3{\leavevmode
  \hbox{\m@@th$\m@th \mathchar"#1#2#3$}}

\def\dots{\relax\ifmmode\ldots\else$\m@th\ldots\,$\relax\fi}

\def\dotfill{\cleaders\hbox{\m@@th$\m@th \mkern1.5mu.\mkern1.5mu$}\hfill}
\def\rightarrowfill{$\m@th\mathord-\mkern-6mu%
  \cleaders\hbox{\m@@th$\mkern-2mu\mathord-\mkern-2mu$}\hfill
  \mkern-6mu\mathord\rightarrow$\relax}
\def\leftarrowfill{$\m@th\mathord\leftarrow\mkern-6mu%
  \cleaders\hbox{\m@@th$\mkern-2mu\mathord-\mkern-2mu$}\hfill
  \mkern-6mu\mathord-$\relax}

\def\downbracefill{$\m@th\braceld\leaders\vrule\hfill\braceru
  \bracelu\leaders\vrule\hfill\bracerd$\relax}
\def\upbracefill{$\m@th\bracelu\leaders\vrule\hfill\bracerd
  \braceld\leaders\vrule\hfill\braceru$\relax}

\def\angle{{\vbox{\m@@th\ialign{$\m@th\scriptstyle##$\crcr
      \not\mathrel{\mkern14mu}\crcr
      \noalign{\nointerlineskip}
      \mkern2.5mu\leaders\hrule height.34pt\hfill\mkern2.5mu\crcr}}}}

\def\big#1{{\m@@th\hbox{$\left#1\vbox to8.5\p@{}\right.\n@space$}}}
\def\Big#1{{\m@@th\hbox{$\left#1\vbox to11.5\p@{}\right.\n@space$}}}
\def\bigg#1{{\m@@th\hbox{$\left#1\vbox to14.5\p@{}\right.\n@space$}}}
\def\Bigg#1{{\m@@th\hbox{$\left#1\vbox to17.5\p@{}\right.\n@space$}}}
\def\n@space{\nulldelimiterspace\z@ \m@th}

\def\root#1\of{\setbox\rootbox\hbox{\m@@th$\m@th\scriptscriptstyle{#1}$}
  \mathpalette\r@@t}
\def\r@@t#1#2{\setbox\z@\hbox{\m@@th$\m@th#1\sqrt{#2}$\relax}
  \dimen@\ht\z@ \advance\dimen@-\dp\z@
  \mkern5mu\raise.6\dimen@\copy\rootbox \mkern-10mu \box\z@}

\def\mathph@nt#1#2{\setbox\z@\hbox{\m@@th$\m@th#1{#2}$}\finph@nt}

\def\mathsm@sh#1#2{\setbox\z@\hbox{\m@@th$\m@th#1{#2}$}\finsm@sh}

\def\@vereq#1#2{\lower.5\p@\vbox{\m@@th\baselineskip\z@skip\lineskip-.5\p@
    \ialign{$\m@th#1\hfil##\hfil$\crcr#2\crcr=\crcr}}}

\def\mathpalette#1#2{\sam\mathchoice{#1\displaystyle{#2}}%
  {#1\textstyle{#2}}{#1\scriptstyle{#2}}{#1\scriptscriptstyle{#2}}\mas}

\def\widehat#1{\setbox\z@\hbox{\sam$#1$}%
 \ifdim\wd\z@>\tw@ em\mathaccent"0\msbfam@5B{#1}%
 \else\mathaccent"0362{#1}\fi}
\def\widetilde#1{\setbox\z@\hbox{\sam$#1$}%
 \ifdim\wd\z@>\tw@ em\mathaccent"0\msbfam@5D{#1}%
 \else\mathaccent"0365{#1}\fi}

 \def\dots{\relax{}
  \ifmmode\def\thedots{\mdots@}\else\def\thedots{\tdots@}\fi %
  \thedots}

 \let\@ldeqno\eqno\let\@ldleqno\leqno
 \def\eqno{\everymath{}\@ldeqno} \def\leqno{\everymath{}\@ldleqno}

  \let\@ldeqalignno\eqalignno
  \def\eqalignno#1{\sam\@ldeqalignno{#1}\mas}
  \let\@ldeqalign\eqalign
  \def\eqalign#1{\sam\@ldeqalign{#1}\mas}

 \def\overrightarrow#1{\vbox{\m@th\ialign{##\crcr
      \rightarrowfill\crcr\noalign{\kern-\p@\nointerlineskip}
      $\hfil\displaystyle{#1}\hfil$\crcr}}}
 \def\overleftarrow#1{\vbox{\m@th\ialign{##\crcr
      \leftarrowfill\crcr\noalign{\kern-\p@\nointerlineskip}
      $\hfil\displaystyle{#1}\hfil$\crcr}}}
 \def\overbrace#1{\mathop{\vbox{\m@th\ialign{##\crcr\noalign{\kern3\p@}
      \downbracefill\crcr\noalign{\kern3\p@\nointerlineskip}
      $\hfil\displaystyle{#1}\hfil$\crcr}}}\limits}
 \def\underbrace#1{\mathop{\vtop{\m@th\ialign{##\crcr
      $\hfil\displaystyle{#1}\hfil$\crcr\noalign{\kern3\p@\nointerlineskip}
      \upbracefill\crcr\noalign{\kern3\p@}}}}\limits}

  \let\@ldmatrix\matrix
  \let\end@ldmatrix\endmatrix
  \def\matrix{\sam\@ldmatrix}
  \def\endmatrix{\end@ldmatrix\mas}
  \let\@ldgather\gather
  \let\end@ldgather\endgather
  \def\gather{\sam\@ldgather}
  \def\endgather{\end@ldgather\mas}
  \let\@ldalign\align
  \let\end@ldalign\endalign
  \def\align{\sam\@ldalign}
  \def\endalign{\end@ldalign\mas}
  \let\@ldaligned\aligned
  \let\end@ldaligned\endaligned
  \def\aligned{\sam\@ldaligned}
  \def\endaligned{\end@ldaligned\mas}
  \let\@ldtag\tag
  \def\tag{\sam\@ldtag}
   %

   \let\MinCDArrowWidth\minCDaw@




\newskip\insertskipamount\newskip\inserthardskipamount
\insertskipamount 6pt plus2pt 
\inserthardskipamount 6pt
\def\insertskip{\vskip\insertskipamount}
\newcount\SplitTest
\def\SetSplitTest{\SplitTest\insertpenalties
  \insert\topins{\floatingpenalty1}%
  \advance\SplitTest-\insertpenalties}
\def\midinsert{\par
 \SaveLastSkip\penalty-150\SetSplitTest\RestoreLastSkip
 \ifnum\SplitTest=-1
  \@midfalse\p@gefalse\else\@midtrue\fi\@ins}
\def\@ins{\par\begingroup\setbox\z@\vbox\bgroup%
  \vglue\inserthardskipamount}
\def\endinsert{\egroup 
  \if@mid \dimen@\ht\z@ \advance\dimen@\dp\z@
    \advance\dimen@\insertskipamount
    \advance\dimen@\pagetotal\advance\dimen@-\pageshrink
    \ifdim\dimen@>\pagegoal\@midfalse\p@gefalse\fi\fi
  \if@mid%
    \ifdim\lastskip<\insertskipamount\removelastskip\insertskip\fi
    \nointerlineskip\box\z@\penalty-200\insertskip
  \else%
    \SaveLastSkip
    \insert\topins{\penalty100 
    \splittopskip\z@skip
    \splitmaxdepth\maxdimen \floatingpenalty\z@
    \ifp@ge \dimen@\dp\z@
    \vbox to\vsize{\unvbox\z@\kern-\dimen@}
    \else \box\z@\nobreak\insertskip\fi}
    \RestoreLastSkip
   \fi\endgroup}


  \newcount\notenumber
  
  \def\note{\advance\notenumber by 1
    \footnote{\the\notenumber)}}

  \newbox\footbox

  \def\footnote#1{\let\@sf\empty
    \ifhmode\edef\@sf{\spacefactor\the\spacefactor}\/\fi
    \sam${}^{\fam0 #1}$\@sf\vfootnote{#1}}%

  \def\vfootnote#1{\insert\footins\bgroup
     \interlinepenalty100 \splittopskip=1pt
     \floatingpenalty=20000
     \leftskip=0pt\rightskip=0pt%
     \parindent=.3em
     \Smallfonts\rm
     \FootItem@{#1}
     \futurelet\next\fo@t}

  \def\FootItem@#1{\par\hangafter1\hangindent=\FootHang
     \setbox0=\hbox{\ignorespaces#1\unskip}%
     \dimen0=.4em\SetOverhang@
     \noindent\rlap{\box0}\kern\Overhang\ignorespaces}


  \def\fo@t{\ifcat\bgroup\noexpand\next \let\next\f@@t
    \else\let\next\f@t\fi \next}
  \def\f@@t{\bgroup\aftergroup\@foot\let\next}
  \def\f@t#1{\baselineskip=10pt\lineskip=1pt
            \lineskiplimit=0pt #1\@foot}%
  \def\@foot{
        \hbox{\vrule height0pt depth5pt width0pt}
        \egroup}
  \skip\footins=12 pt plus 0pt minus 0pt 
  \count\footins=1000 
  \dimen\footins=8in 



 \def\osumess#1{\EdSpider{\immediate\write16{Line \the\inputlineno: #1}}}%
 \def\HideEdStuff{\gdef\EdSpider##1{}}

 \font\BigSym=cmmi10 scaled \magstep 4

 \def\change{\InLMargin{\hbox{\BigSym \char63\kern10pt}}}

 \def\beginchange{\InLMargin{\hbox{\sam\twelvepoint$\heartsuit$\kern10pt}}}

 \def\endchange{\InLMargin{\hbox{\sam\twelvepoint$\spadesuit$\kern10pt}}}

 \def\InLMargin#1{\strut\vadjust{%
     \kern-\strutdepth
     \vtop to \strutdepth{%
         \baselineskip\strutdepth
         \llap{\sam$\smash{\hbox{\EdSpider{#1}}}$}\null}}}

 \def\strutdepth{\dp\strutbox}
 \def\strutheight{\ht\strutbox}

 \def\NoteInRMargin#1{\strut\vadjust{%
     \kern-1.001\strutdepth
     \vtop to \strutdepth{%
       \baselineskip\strutdepth
       \vss\rlap{\ninepoint\unskip\hskip\hsize
         \vtop to 0pt{%
           \hsize=16em\hfuzz=\hsize
           \leftskip=10pt%
           \rightskip=0pt plus 10000pt%
           \baselineskip=9.8pt\lineskip=.2pt%
           \let\\\break
           \noindent\EdSpider{#1}\vss}%
                \kern10pt}\hbox{}}
       }}

 \def\ednote#1{\NoteInRMargin{\tentt #1}}

 \def\cbar{\InLMargin{%
      \dimen0=\strutdepth\advance\dimen0 by \lineskip
      \vrule width 3pt
      height \strutheight depth \dimen0 \kern
      3pt}}

 \def\ccbar{\InLMargin{%
      \dimen0=2\strutdepth\advance\dimen0 by 2\lineskip
      \vrule width 3pt
        height 3\strutheight depth \dimen0 \kern
      3pt}}

 \newinsert\TRMargIns
 \dimen\TRMargIns=\maxdimen

  \def\Ednote#1{\insert\TRMargIns{%
       \vbox to 0pt{\hsize=140pt\hfuzz=\hsize
           \leftskip=6pt%
           \rightskip=0pt plus 10000pt%
           \baselineskip=9.8pt\lineskip=.2pt%
           \let\\\break
           \SetPageRemainder
           \vglue540pt\vglue-\PageRemainder
           \noindent\EdSpider{\tentt #1}\vss}%
       \smallskip}}

 \def\KillEdStuff{\def\ednote##1{}\def\Ednote##1{}%
      \let\change\relax\let\beginchange\relax\let\endchange\relax
       \let\cbar\relax\let\ccbar\relax}


  \topskip=12pt
  \newskip\StdBaselineskip 
  \StdBaselineskip 12pt
  \lineskip=1.1pt
  \lineskiplimit=.8pt
  \widowpenalty=10000 
  \clubpenalty=10000  
  \abovedisplayskip=6pt plus 1pt minus 1pt
  \abovedisplayshortskip=3pt plus 1.5pt
  \belowdisplayskip=6pt plus 1pt minus 1pt
  \belowdisplayshortskip=5pt plus 1pt minus 1pt
  \hfuzz=1.5pt   

  \def\StdPretolerance{100}
  \tolerance=\StdPretolerance

  \newdimen\StdMathsurround
  \StdMathsurround=1.5pt 
  \mathsurround=\StdMathsurround
  \Mas                   

   \def\prose{\relax\hbox{\kern.6\StdMathsurround}}
  
  \def\StdParskip{0pt}    
  \parskip=\StdParskip
  \parindent=0.5cm
 

  \def\Times{ptmr  } 
  \def\TimesI{ptmri  } 
  \def\TimesB{ptmb  }
  \def\TimesBI{ptmbi  }
  \def\HelveticaN{phvrrn }

  =\Times at 10bp
  =\TimesB at 10bp
  \font\tenit=\TimesI at 10bp
  =\TimesBI at 10bp

  \font\tenmrm=cmr10  


    =\Times at 9bp 
    \font\nineit=\TimesI at 9bp 
    =\TimesB at 9bp 
    =\TimesBI at 9bp 

    =\HelveticaN at 9bp 


  =\Times at 12bp
  \font\twelveit=\TimesI at 12bp
  =\TimesB at 12bp


  \font\titleit=\TimesI at 14.4bp
  =\TimesB at 14.4bp

 \SetAuthorHead{AuthorHead} 
 \SetTitleHead{TitleHead}  


  \def\lBr{\raise.125ex\hbox{[\kern.1125ex}}
  \def\rBr{\raise.125ex\hbox{\kern.1125ex]}}

 \setbox\footbox=\hbox{\Smallfonts 2)~}



  \bgroup
  \catcode`\@=11 
  \gdef\itSpacing{%
     \xspaceskip=.31em plus.1em minus.05em \sfcode `f=2001
     \itWarning@\let\itWarning@\itWarning@@}
  \gdef\itSpacingOff{%
     \xspaceskip=0pt \sfcode `f=1000
     \let\itWarning@\relax}
   \global\let\itWarning@\relax
  \gdef\itWarning@@{\errmessage{%
  Special italic spacing already in force
  (you have probably omitted an ``endth'').
  See itSpacing macro in osuPSfnt.sty
         }}
  \egroup

 \fontdimen1\titlebf=0.0pt
 \fontdimen2\titlebf=3.6135pt
 \fontdimen3\titlebf=2.8908pt
 \fontdimen4\titlebf=1.44539pt
 \fontdimen5\titlebf=6.64882pt
 \fontdimen6\titlebf=14.45398pt
 \fontdimen7\titlebf=1.60439pt

 \fontdimen1\tenbi=0.26794pt
 \fontdimen2\tenbi=2.50937pt
 \fontdimen3\tenbi=2.00749pt
 \fontdimen4\tenbi=1.00374pt
 \fontdimen5\tenbi=4.59717pt
 \fontdimen6\tenbi=10.03749pt
 \fontdimen7\tenbi=1.11415pt

 \fontdimen1\twelverm=0.0pt
 \fontdimen2\twelverm=3.01125pt
 \fontdimen3\twelverm=2.409pt
 \fontdimen4\twelverm=1.2045pt
 \fontdimen5\twelverm=5.39615pt
 \fontdimen6\twelverm=12.045pt
 \fontdimen7\twelverm=1.33699pt

 \fontdimen1\twelveit=0.27731pt
 \fontdimen2\twelveit=3.01125pt
 \fontdimen3\twelveit=2.409pt
 \fontdimen4\twelveit=1.2045pt
 \fontdimen5\twelveit=5.37207pt
 \fontdimen6\twelveit=12.045pt
 \fontdimen7\twelveit=1.33699pt

 \fontdimen1\twelvebf=0.0pt
 \fontdimen2\twelvebf=3.01125pt
 \fontdimen3\twelvebf=2.409pt
 \fontdimen4\twelvebf=1.2045pt
 \fontdimen5\twelvebf=5.5407pt
 \fontdimen6\twelvebf=12.045pt
 \fontdimen7\twelvebf=1.33699pt

 \fontdimen1\tenrm=0.0pt
 \fontdimen2\tenrm=2.50937pt
 \fontdimen3\tenrm=2.00749pt
 \fontdimen4\tenrm=1.00374pt
 \fontdimen5\tenrm=4.49678pt
 \fontdimen6\tenrm=10.03749pt
 \fontdimen7\tenrm=1.11415pt

 \fontdimen1\tenit=0.27731pt
 \fontdimen2\tenit=2.50937pt
 \fontdimen3\tenit=2.00749pt
 \fontdimen4\tenit=1.00374pt
 \fontdimen5\tenit=4.47672pt
 \fontdimen6\tenit=10.03749pt
 \fontdimen7\tenit=1.11415pt

 \fontdimen1\tenbf=0.0pt
 \fontdimen2\tenbf=2.50937pt
 \fontdimen3\tenbf=2.00749pt
 \fontdimen4\tenbf=1.00374pt
 \fontdimen5\tenbf=4.61723pt
 \fontdimen6\tenbf=10.03749pt
 \fontdimen7\tenbf=1.11415pt

 \fontdimen1\ninerm=0.0pt
 \fontdimen2\ninerm=2.25842pt
 \fontdimen3\ninerm=1.80673pt
 \fontdimen4\ninerm=0.90337pt
 \fontdimen5\ninerm=4.0471pt
 \fontdimen6\ninerm=9.03374pt
 \fontdimen7\ninerm=1.00273pt

 \fontdimen1\nineit=0.27731pt
 \fontdimen2\nineit=2.25842pt
 \fontdimen3\nineit=1.80673pt
 \fontdimen4\nineit=0.90337pt
 \fontdimen5\nineit=4.02904pt
 \fontdimen6\nineit=9.03374pt
 \fontdimen7\nineit=1.00273pt

 \fontdimen1\ninebf=0.0pt
 \fontdimen2\ninebf=2.25842pt
 \fontdimen3\ninebf=1.80673pt
 \fontdimen4\ninebf=0.90337pt
 \fontdimen5\ninebf=4.15552pt
 \fontdimen6\ninebf=9.03374pt
 \fontdimen7\ninebf=1.00273pt


 \newcount\MaxSpaceFactor
 \MaxSpaceFactor=3000 

 \def\ItemStyle{\rm}
 \def\NrStyle{\rm}
 \def\ItemItemStyle{\rm}

 \MaxItemTag{(iii)}
 \MaxItemItemTag{(iii)}
 \MaxNrTag{(2)}
 \MaxFootTag{2)}
 \def\ReferenceHang{30pt}

 \catcode`\@=\active


\loadbold

=\Times  
=\Times scaled750
=\Times scaled650
\font\rms=\Times scaled 920 

=\TimesBI scaled 860
=\TimesI scaled 860

\textfont0=\rrm  
\scriptfont0=\erm 
\scriptscriptfont0=\srm

\def\Augment#1#2{%
    \toks0\expandafter{#1}\toks2{#2}%
    \edef#1{\the\toks0\the\toks2}}

 \font\twelverma=\Times  scaled 1200
 \font\tenrma=\Times  scaled 1000
 \font\ninerma=\Times scaled 920
 =\Times scaled 840
 \font\sevenrma=\Times scaled 760
 =\Times scaled 680
 \font\fiverma=\Times scaled 600

 \Augment\tenpoint{%
  \textfont0=\tenrma  \scriptfont0=\sevenrma  
  \scriptscriptfont0=\fiverma  }

 \Augment\ninepoint{%
  \textfont0=\ninerma  \scriptfont0=\sevenrma 
  \scriptscriptfont0=\fiverma}

 \Augment\twelvepoint{%
  \textfont0=\twelverma  \scriptfont0=\ninerma  
  \scriptscriptfont0=\sevenrma}

\mathsurround=1pt
\hsize=13.45truecm
\vsize=19.5truecm
\hoffset=1.25truecm
\voffset=2truecm
\advance\baselineskip by 2pt

\predefine\til{\~}
\def\~#1{\relax\ifmmode\widetilde{#1}\else\til{#1}\fi}

\redefine \le{\leqslant}

\define \wt#1{\mathaccent"0365{#1}}
\define \wh#1{\mathaccent"0362{#1}}

\define \iss{\,\Mathaccent{\raise -.8 ex\hbox{$\widetilde{}$\kern.1em}}\rightarrow\,}

\define \inlim{{\varinjlim}\vphantom{i}\,}

\define \ur{\mathop{\fam0 ur}}
\define \nr{\mathop{\fam0 ur}}

\define \ab{\mathop{\fam0 ab}}

\define \sep{\mathop{\fam0 sep}}

\define \chr{\mathop{\fam0 char}\,}

\define \res{\operatorname{\fam0 res}}

\define \Br{\operatorname{\fam0 Br}}

\define \cor{\operatorname{\fam0 cor}}

\define \Gal{\mathop{\fam0 Gal}}
\define \Hom{\operatorname{\fam0 Hom}}

\define \inv{\mathop{\fam0 inv}}

\define \gr{\operatorname{\fam0 gr}\!}

\Mas
\HideEdStuff
\rm 
 

\def\issn{{\nineit ISSN 1464-8997 (on line) 1464-8989 (printed)}}

\def\gtp{{\nineit Published 10 December 2000: \ \copyright\ Geometry \& 
Topology Publications}}

\def\gtv3{{\nineit Geometry \& Topology Monographs, Volume 3 (2000) --
Invitation to higher local fields}}


\def\lione
{{\rms Geometry \& Topology Monographs}}

\def \litwo{{\rms Volume 3: Invitation to higher local fields
}} 

\def\tinfo #1.#2.#3-#4
{{
\noindent  {\lione} \hfill 
\par 
\vskip-1.5pt
\noindent {\litwo} \hfill
\par 
\vskip-1,5pt
\noindent {\rms Part #1, section #2, pages #3--#4} \hfill
\vskip24pt 
}}

\def\tinfos #1.#2.#3-#4
{{
\noindent  {\lione} \hfill 
\par 
\vskip-1.5pt
\noindent {\litwo} \hfill
\par 
\vskip-1.5pt
\noindent {\rms Pages #3--#4} \hfill
\vskip24pt 
}}

\def\tinfoi #1
{{
\noindent  {\lione} \hfill 
\par 
\vskip-1.5pt
\noindent {\litwo} \hfill
\par 
\vskip-1.5pt
\noindent {\rms Pages iii--xi: Introduction and contents} \hfill
\vskip26pt 
}}


  \def\titlepagehead{\hfil}

  \newif\iftitlepage\titlepagefalse
  \newif\ifblankpage\blankpagefalse
  \def\makeheadline{
     \ifblankpage{}\else%
     \iftitlepage
\vbox{\line{\vbox to 8.5pt{}
\ninerm
\copy\HLinebox \hfill
\hglue5mm\ninebf\folio 
\titlepagehead}}%
      \else
\vbox{\ifodd\pageno\rightheadline\else\leftheadline\fi}%
      \fi\vskip 12pt\fi}%
     \def\rightheadline{\line{\vbox to 8.5pt{}%
      \ninerm
\copy\TitleBox \hfill
\hglue5mm\ninebf\folio}}%
     \def\leftheadline{\line{\vbox to 8.5pt{}%
        \unskip\ninerm\unskip\ninebf\folio\hglue5mm
 \hfill \copy\AuthorBox
}}

 \footline={\ifblankpage{}\else
\iftitlepage\ninepoint\sam\hfill
\line{\vbox to 8.5pt{}
\copy\TFLinebox
\hfill
\hglue5mm 
}
            \else
\ninepoint\sam\hfill
\line{\vbox to 8.5pt{}
\copy\FLinebox
\hfill 
\hglue5mm
}
\hfil\fi\global\titlepagefalse\fi}

\def\blankpage{{\blankpagetrue\noindent\hbox to 10pt{\hss}\vfill
\pagebreak}}

\tenpoint\rm 
 

\pageno=53

\tinfo I.5.53-60

\SetTFLinebox{\gtp }
\SetFLinebox{\gtv3 }
\SetHLinebox{\issn}

\H 5. Kato's higher local class field theory

Masato Kurihara

\SetAuthorHead{M. Kurihara}
\SetTitleHead{Part I. Section 5. Kato's higher local class field
theory\qquad\qquad}

\HH 5.0. Introduction

We first recall the classical local class field theory. 
Let $K$ be a finite extension of ${\Bbb Q}_{p}$ or ${\Bbb F}_{q}((X))$. 
The main theorem of local class field theory consists of 
the isomorphism theorem and existence theorem. 
In this section we consider the isomorphism theorem. 

An outline of one of the proofs is as follows.
First, for the Brauer group $\Br(K)$, 
an isomorphism 
$$\inv\colon \Br(K) \iss  {\Bbb Q}/{\Bbb Z}$$
is established; 
it mainly follows from an isomorphism 
$$H^{1}(F, {\Bbb Q}/{\Bbb Z}) 
\iss  {\Bbb Q}/{\Bbb Z}$$
where $F$ is the residue field of $K$. 

Secondly, we denote by $X_{K}=\Hom_{\text{\rm cont}}(G_{K}, {\Bbb Q}/{\Bbb
Z})$ 
the group of 
continuous homomorphisms 
from $G_{K}=\Gal(\Overline{K}/K)$ to ${\Bbb Q}/{\Bbb Z}$. 
We consider a pairing 
$$K^{{*}} \times X_{K} \longrightarrow {\Bbb Q}/{\Bbb Z}$$
$$(a, \chi) \mapsto \inv(\chi,a)$$
where $(\chi,a)$ is the cyclic algebra associated with $\chi$ and $a$. 
This pairing induces a homomorphism 
$$\Psi_{K}\colon K^{{*}} \longrightarrow \Gal(K^{\ab}/K)=
\Hom(X_{K},{\Bbb Q}/{\Bbb Z})$$ 
which is called the reciprocity map. 

Thirdly, for a finite abelian extension $L/K$, 
we have a diagram 
$$
\CD
L^{{*}} @>{\Psi_{L}}>> \Gal(L^{\ab}/L) \\ 
@V N VV @VVV \\
K^{{*}}  @>{\Psi_{K}}>>
 \Gal(K^{\ab}/K) 
\endCD
$$
which is commutative by 
the definition of the reciprocity maps. 
Here, $N$ is the norm map and the right vertical map is the canonical map. 
This induces a homomorphism 
$$\Psi_{L/K}\colon K^{{*}}/NL^{{*}} \longrightarrow \Gal(L/K).$$
The isomorphism theorem tells us that the above map 
is bijective.  

To show the bijectivity of $\Psi_{L/K}$, we can reduce to the case 
where $|L:K|$ is a prime $\ell$. 
In this case, the bijectivity follows immediately from 
a famous exact sequence 
$$L^{{*}} @>{N}>> K^{{*}} 
@>{\cup \chi}>> 
\Br(K) @>{\res}>> \Br(L)$$ 
for a cyclic extension $L/K$ 
(where $\cup \chi$ is the cup product with $\chi$, and 
$\res$ is the restriction map). 

In this section  we sketch a proof of the isomorphism theorem
for 
a higher dimensional local field as an analogue of the above argument. 
For the existence theorem see the paper by Kato in this volume
and subsection 10.5. 

\HH 5.1. Definition of $H^{q}(k)$

In this subsection, for any field $k$ and $q > 0$, we recall the
definition of the cohomology 
group $H^{q}(k)$ (\cite{K2}, see also subsections  2.1 and 2.2
and A1 in the appendix to section~2). 
If $\chr(k)=0$, 
we define $H^{q}(k)$ as a Galois cohomology group 
$$H^{q}(k)=H^{q}(k, {\Bbb Q}/{\Bbb Z}(q-1))$$
where $(q-1)$ is the $(q-1)$st Tate twist. 

If $\chr(k)=p>0$, 
then following Illusie \cite{I}  
we define 
$$H^{q}(k, {\Bbb Z}/p^{n}(q-1))=
H^{1}(k, W_{n}\Omega_{k^{\sep}, \log}^{q-1}).$$ 
We can explicitly describe 
$H^{q}(k, {\Bbb Z}/p^{n}(q-1))$ as the 
group isomorphic to 
$$W_{n}(k) \otimes (k^{{*}})^{\otimes (q-1)}/J$$ 
where $W_{n}(k)$ is the ring of Witt vectors of length $n$, and 
$J$ is the subgroup generated by  elements of the form 
\Roster
\Item{} $w \otimes b_{1} \otimes \dots \otimes b_{q-1}$ such that
$b_{i}=b_{j}$ for 
some $i \neq j$, and 
\Item{} $(0,\dots ,0,a,0,\dots,0) \otimes a \otimes b_1 \otimes \cdots
\otimes b_{q-2}$,
and 
\Item{} $({\bold F}-1)(w) \otimes b_{1} \otimes \cdots \otimes b_{q-1}$ 
($\bold F$ is the Frobenius map on Witt vectors). 
\endRoster  

We define 
$H^{q}(k, {\Bbb Q}_{p}/{\Bbb Z}_{p}(q-1))=
\inlim H^{q}(k, {\Bbb Z}/p^{n}(q-1))$, and define 
$$H^{q}(k)= \bigoplus_{\ell} 
H^{q}(k, {\Bbb Q}_{\ell}/{\Bbb Z}_{\ell}\,(q-1))$$
where $\ell$ ranges over all prime numbers. 
(For $\ell \neq p$, the right hand side is the usual Galois cohomology 
of the $(q-1)$st Tate twist of ${\Bbb Q}_{\ell}/{\Bbb Z}_{\ell}$.) 

Then for any $k$ we have 
$$
\alignat2
& H^{1}(k)=X_{k} \qquad &&\text{($X_{k}$ is as in 5.0, 
the group of characters)}, \\
& H^{2}(k)=\Br(k) \qquad && \text{(Brauer group)}.
\endalignat 
$$

We explain the second equality in  the case of $\chr(k)=p>0$. 
The relation between the Galois cohomology group and the Brauer group 
is well known, so we consider only the $p$-part. 
By our definition, 
$$H^{2}(k, {\Bbb Z}/p^{n}(1))=
H^{1}(k, W_{n}\Omega_{k^{\sep}, \log}^{1}).$$ 
From the bijectivity of the differential symbol (Bloch--Gabber--Kato's
theorem 
in subsection A2 in the appendix to section~2), we have 
$$H^{2}(k, {\Bbb Z}/p^{n}(1))=
H^{1}(k, (k^{\sep})^{{*}}/((k^{\sep})^{{*}})^{p^{n}}).$$ 
From the exact sequence 
$$0 \longrightarrow (k^{\sep})^{{*}} @>{p^{n}}>> 
(k^{\sep})^{{*}} \longrightarrow 
(k^{\sep})^{{*}}/((k^{\sep})^{{*}})^{p^{n}} \longrightarrow 0$$ 
and an isomorphism $\Br(k)=H^{2}(k, (k^{\sep})^{{*}})$, 
$H^{2}(k, {\Bbb Z}/p^{n}(1))$ is isomorphic to the $p^{n}$-torsion 
points of $\Br(k)$. 
Thus, we get $H^{2}(k)=\Br(k)$. 

\medskip

If $K$ is a henselian discrete valuation field 
with residue field $F$, 
we have a canonical map 
$$i_{F}^{K}\colon H^{q}(F) \longrightarrow H^{q}(K).$$
If $\chr(K)=\chr(F)$, this map is defined naturally from the definition 
of $H^{q}$ 
(for the Galois cohomology part, we use a natural map 
$\Gal(K^{\sep}/K) \longrightarrow \Gal(K_{\ur}/K)=\Gal(F^{\sep}/F)$) . 
If $K$ is of mixed characteristics $(0,p)$, the prime-to-$p$-part is
defined 
naturally and the $p$-part is 
defined as follows. 
For the class $[w \otimes b_{1} \otimes \dots \otimes b_{q-1}]$ in 
$H^{q}(F,{\Bbb Z}/p^{n}(q-1))$ we define 
$i_{F}^{K}([w \otimes b_{1} \otimes \dots \otimes b_{q-1}])$ as 
the class of 
$$i(w) \otimes \wt{b_{1}} \otimes \dots \otimes \wt{b_{q-1}}$$ 
in 
$H^{1}(K,{\Bbb Z}/p^{n}(q-1))$, where 
$i\colon W_{n}(F) \rightarrow H^{1}(F,{\Bbb Z}/p^{n}) \rightarrow 
H^{1}(K,{\Bbb Z}/p^{n})$ is the composite of the map 
given by Artin--Schreier--Witt theory 
and the canonical map, and $\wt{b_{i}}$ is a lifting of $b_{i}$ to 
$K$. 

\th Theorem {{\rm(Kato \cite{K2, Th. 3})}}

Let $K$ be a henselian discrete valuation field, $\pi$ be a prime element, 
and $F$ be the residue field. 
We consider a homomorphism 
$$\aligned
&i=(i_{F}^{K}, i_{F}^{K} \cup \pi) \colon  
H^{q}(F) \oplus H^{q-1}(F) \longrightarrow H^{q}(K)\\
&(a,b) \mapsto i_{F}^{K}(a) + i_{F}^{K}(b) \cup \pi
\endaligned
$$
where $i_{F}^{K}(b) \cup \pi$ is the element obtained from the pairing 
$$H^{q-1}(K) \times K^{{*}} \longrightarrow H^{q}(K)$$
which is defined by Kummer theory and the cup product, and 
the explicit description of $H^{q}(K)$ in the  case of $\chr(K) >0$. 
Suppose $\chr(F)=p$. 
Then $i$ is bijective in the prime-to-$p$ component. 
In the $p$-component, $i$ is injective and the image coincides with the 
$p$-component of the kernel of $H^{q}(K) \rightarrow H^{q}(K_{\nr})$ 
where $K_{\nr}$ is the maximal unramified extension of $K$. 
\endth

From this theorem and Bloch--Kato's theorem in section~4, we obtain 

\th Corollary

Assume that $\chr(F)=p>0$, $|F:F^{p}|=p^{d-1}$, and that there is an 
isomorphism 
$H^{d}(F)\iss  {\Bbb Q}/{\Bbb Z}$. 

Then, $i$ induces an isomorphism 
$$H^{d+1}(K)\iss  {\Bbb Q}/{\Bbb Z}.$$
\endth

A typical example which satisfies the assumptions of the above corollary 
is a $d$-di\-me\-n\-si\-o\-nal local field 
(if the last residue field is quasi-finite (not necessarily finite), 
the assumptions are satisfied).

\HH 5.2. Higher dimensional local fields

We assume that $K$ is a $d$-dimensional local field, and 
$F$ is the residue field of $K$, which is a $(d-1)$-dimensional local
field. 
Then, by the corollary in the previous subsection and induction on $d$, 
there is a canonical isomorphism 
$$\inv\colon H^{d+1}(K)\iss  {\Bbb Q}/{\Bbb Z}.$$
This corresponds to the first step of the proof of the classical 
isomorphism theorem which we described in the introduction. 

The cup product defines a pairing 
$$K_{d}(K) \times H^{1}(K) \longrightarrow H^{d+1}(K) 
\simeq {\Bbb Q}/{\Bbb Z}.$$
This pairing induces a homomorphism 
$$\Psi_{K}\colon K_{d}(K) \longrightarrow \Gal(K^{\ab}/K)\simeq 
\Hom(H^{1}(K),{\Bbb Q}/{\Bbb Z})$$ 
which we call the {\it reciprocity map}. 
Since the isomorphism 
$\inv\colon H^{d}(K) \longrightarrow {\Bbb Q}/{\Bbb Z}$ is 
naturally constructed, 
for a finite abelian extension $L/K$  
we have a commutative diagram 
$$
\CD
H^{d+1}(L) @>{\inv_{L}}>> {\Bbb Q}/{\Bbb Z} \\ 
@V{\cor}VV  @VVV \\
H^{d+1}(K)  @>{\inv_{K}}>> {\Bbb Q}/{\Bbb Z}. 
\endCD
$$
So the diagram 
$$
\CD
K_{d}(L) @>{\Psi_{L}}>> \Gal(L^{\ab}/L) \\ 
@V{N}VV  @VVV \\
K_{d}(K)  @>{\Psi_{K}}>> 
 \Gal(K^{\ab}/K) 
\endCD
$$
is commutative 
where $N$ is the norm map and the right vertical map is the canonical map. 
So, as in the classical case, we have a homomorphism 
$$\Psi_{L/K}\colon K_{d}(K)/NK_{d}(L) \longrightarrow \Gal(L/K).$$

\th Isomorphism Theorem 

$\Psi_{L/K}$ is an isomorphism. 
\endth

We outline a proof. 
We may assume that 
$L/K$ is cyclic of degree $\ell$. 
As in the classical case in the introduction, we may study a 
sequence 
$$K_{d}(L) @>{N}>> K_{d}(K) 
@>{\cup \chi}>> H^{d+1}(K) 
@>{\res}>> H^{d+1}(L),$$
but here we describe a more elementary proof. 

\smallskip

First of all, 
using the argument in \cite{S, Ch.5} 
by calculation of symbols 
one can obtain 
$$|K_{d}(K):NK_{d}(L)| \le \ell.$$ 
We outline a proof of this inequality. 

It is easy to see that it is sufficient to consider
the case of prime $\ell$.
(For another calculation of the index of the norm group see subsection~6.7).

Recall that $K_{d}(K)$ has a filtration 
$U_{m}K_{d}(K)$ 
as in subsection 4.2. 
We consider 
$\gr_{m}K_{d}(K)=U_{m}K_{d}(K)/U_{m+1}K_{d}(K)$. 

If $L/K$ is unramified, the norm map 
$N\colon K_{d}(L) \rightarrow K_{d}(K)$ induces 
surjective homomorphisms 
$\gr_{m}K_{d}(L) \rightarrow \gr_{m}K_{d}(K)$ 
for all $m >0$. 
So $U_{1}K_{d}(K)$ is in $NK_{d}(L)$. 
If we denote by $F_{L}$ and $F$ the residue fields of 
$L$ and $K$ respectively, 
the norm map induces a surjective homomorphism 
$K_{d}(F_{L})/\ell \rightarrow K_{d}(F)/\ell$ 
because $K_{d}(F)/\ell$ is isomorphic to 
$H^{d}(F, {\Bbb Z}/\ell(d))$ (cf. sections 2 and 3) 
and the cohomological 
dimension of $F$ \cite{K2, p.220} is $d$. 
Since $\gr_{0}K_{d}(K)=K_{d}(F) \oplus K_{d-1}(F)$ 
(see subsection~4.2), 
the above implies that 
$K_{d}(K)/NK_{d}(L)$ is isomorphic to 
$K_{d-1}(F)/NK_{d-1}(F_{L})$, which is 
isomorphic to $\Gal(F_{L}/F)$ by  class field theory 
of $F$ (we use induction on $d$). 
Therefore $|K_{d}(K):NK_{d}(L)| = \ell$. 

If $L/K$ is totally ramified and $\ell$ is prime to $\chr(F)$, 
by the same argument (cf. the argument in \cite{S, Ch.5}) as above, 
we have $U_{1}K_{d}(K) \subset NK_{d}(L)$. 
Let $\pi_{L}$ be a prime element of $L$, and 
$\pi_{K}=N_{L/K}(\pi_{L})$. 
Then the element $\{\alpha_1,...,\alpha_{d-1},\pi_{K}\}$ 
for $\alpha_i\in K^*$ is in 
$NK_{d}(L)$, so 
$K_{d}(K)/NK_{d}(L)$ is isomorphic to 
$K_{d}(F)/\ell$, which is 
isomorphic to $H^{d}(F, {\Bbb Z}/\ell(d))$, 
so the order is $\ell$. 
Thus, in this case we also have 
$|K_{d}(K):NK_{d}(L)| = \ell$. 

Hence, we may assume $L/K$ is not unramified and is of degree $\ell=p=\chr(F)$. 
Note that $K_{d}(F)$ is $p$-divisible because 
of $\Omega_{F}^{d}=0$ and the bijectivity of the 
differential symbol. 

Assume that $L/K$ is totally ramified. 
Let $\pi_{L}$ be a prime element of $L$, and 
$\sigma$ a generator of $\Gal(L/K)$, and put 
$a=\sigma(\pi_{L})\pi_{L}^{-1}-1$, 
$b=N_{L/K}(a)$, 
and 
$v_{K}(b-1)=i$. 
We study the induced maps 
$\gr_{\psi(m)}K_{d}(L) \rightarrow \gr_{m}K_{d}(K)$ 
from the norm map $N$ on the subquotients 
by the argument in \cite{S, Ch.5}. 
We have 
$U_{i+1}K_{d}(K) \subset NK_{d}(L)$, and 
can show that there is a surjective homomorphism 
(cf. 
\cite{K1, p.669}) 
$$\Omega_{F}^{d-1} \longrightarrow K_{d}(K)/NK_{d}(L)$$ 
such that 
$$x d \log y_{1} \wedge ... \wedge d \log y_{d-1} 
\mapsto \{1+\widetilde{x}b,\widetilde{y_1},...,\widetilde{y_{d-1}}\}$$ 
($\widetilde{x}, \widetilde{y_i}$ are liftings of $x$ and $y_{i}$). 
Furthermore, from 
$$N_{L/K}(1+xa) \equiv 1+(x^{p}-x)b \pmod{U_{i+1}K^{{*}}},$$ 
the above map 
induces 
a surjective homomorphism 
$$\Omega_{F}^{d-1}/
\bigl(({\bf F}-1)\Omega_{F}^{d-1} +d \Omega_{F}^{d-2}\bigr)
\longrightarrow K_{d}(K)/NK_{d}(L).$$
The source group is isomorphic to $H^{d}(F, {\Bbb Z}/p(d-1))$ which is 
of order $p$. 
So we obtain 
$|K_{d}(K):NK_{d}(L)| \le p$. 

Now assume that  $L/K$ is ferociously ramified, i.e.  
$F_{L}/F$ is purely inseparable of degree $p$. We can
use an  argument  similar to the previous one.  
Let $h$ be an element of ${\Cal O}_{L}$ such that $F_{L}=F(\overline{h})$ 
($\overline{h}=h \mod \Cal M_{L}$). 
Let $\sigma$ be a generator of $\Gal(L/K)$, and 
put $a=\sigma(h)h^{-1}-1$, and $b=N_{L/K}(a)$. 
Then we have a surjective homomorphism
(cf. \cite{K1, p.669})  
$$\Omega_{F}^{d-1}/\bigl(({\bf F}-1)\Omega_{F}^{d-1} + d \Omega_{F}^{d-2}\bigr)
\longrightarrow K_{d}(K)/NK_{d}(L)$$ 
such that $$x d \log y_{1} \wedge ... \wedge d \log y_{d-2} \wedge 
d \log N_{F_{L}/F}(\overline{h})  
\mapsto \{1+\widetilde{x}b, \widetilde{y_1},..., \widetilde{y_{d-2}}, \pi \}$$ 
($\pi$ is a prime element of $K$).  
So we get 
$|K_{d}(K):NK_{d}(L)| \le p$. 

\medskip

So in order to obtain the bijectivity of $\Psi_{L/K}$, we have only to 
check the surjectivity. 
We consider the most interesting case $\chr(K)=0$, $\chr(F)=p>0$, and 
$\ell=p$. 
To show the surjectivity of $\Psi_{L/K}$, we have to show 
that there is an element $x \in K_{d}(K)$ such that $\chi \cup x \neq 0$ 
in $H^{d+1}(K)$ where $\chi$ is a character corresponding to $L/K$.
We may assume a primitive $p$-th root of unity is in $K$. 
Suppose that $L$ is given by an equation $X^{p}=a$ for some 
$a \in K\setminus K^p$. 
By Bloch--Kato's theorem (bijectivity of the cohomological symbols in 
section~4), we identify the kernel of  multiplication by $p$ 
on $H^{d+1}(K)$ with 
$H^{d+1}(K,{\Bbb Z}/p(d))$, and with 
$K_{d+1}(K)/p$. 
Then  
our aim is to show that 
there is an element 
$x \in K_{d}(K)$ such that $\{x,a\} \neq 0$ in 
$k_{d+1}(K)=K_{d+1}(K)/p$. 
(Remark. The pairing $K_1(K)/p\times K_d(K)/p
\to K_{d+1}(K)/p$ coincides up to a sign
with Vostokov's symbol defined in subsection~8.3 and 
the latter is non-degenerate which
 provides an alternative proof).

We use the notation of section 4. 
By the Proposition in subsection 4.2, 
we have $$K_{d+1}(K)/p=k_{d+1}(K)=U_{e'}k_{d+1}(K)$$ 
where $e'=v_{K}(p)p/(p-1)$. 
Furthermore, by the same proposition  
there is an isomorphism 
$$H^{d}(F,{\Bbb Z}/p(d-1))=\Omega_{F}^{d-1}/\bigl(({\bf F}-1)\Omega_{F}^{d-1} 
+ d \Omega_{F}^{d-2}\bigr) \longrightarrow k_{d+1}(K)$$ 
such that 
$$x d \log y_{1} \wedge ... \wedge d \log y_{d-1} 
\mapsto \{1+\widetilde{x}b,\widetilde{y_1},...,\widetilde{y_{d-1}}, \pi \}$$ 
where $\pi$ is a uniformizer, and 
$b$ is a certain element of $K$ such that $v_{K}(b)=e'$. 
Note that $H^{d}(F, {\Bbb Z}/p(d-1))$ is of order $p$. 
 
This shows that 
for any uniformizer $\pi$ of $K$, and for any 
lifting $t_{1}$,...,$t_{d-1}$ of a $p$-base of $F$, 
there is an element $x \in {\Cal O}_{K}$ such that  
$$\{1+\pi^{e'}x, t_{1},...,t_{d-1}, \pi\} \neq 0$$ 
in $k_{d+1}(K)$. 

If the class of $a$ is not in 
$U_{1}k_{1}(K)$, 
we may assume $a$ is a uniformizer or 
$a$ is a part of a lifting of a $p$-base of $F$. 
So it is easy to see by the above property that 
there exists an $x$ such that $\{a,x\} \neq 0$. 
If the class of $a$ is in 
$U_{e'}k_{1}(K)$, it is also easily seen 
from the description of 
$U_{e'}k_{d+1}(K)$ that 
there exists an $x$ such that $\{a,x\} \neq 0$.

Suppose $a \in U_{i}k_{1}(K) \setminus U_{i+1}k_{1}(K)$ such that 
$0 < i < e'$. 
We write $a=1+\pi^{i} a'$ for a prime element $\pi$ and 
$a' \in {\Cal O}_{K}^{{*}}$. 
First, we assume that $p$ does not divide $i$. 
We use a formula (which holds in $K_{2}(K)$) 
$$\{1-\alpha,1-\beta\}=\{1-\alpha \beta,-\alpha\}+\{1-\alpha
\beta,1-\beta\}
-\{1-\alpha \beta,1-\alpha\}$$ 
for $\alpha \neq 0,1$, and $\beta \neq 1,\alpha^{-1}$. 
From this formula we have in $k_{2}(K)$ 
$$\{1+\pi^{i}a', 1+\pi^{e'-i}b \}= \{1+ \pi^{e'}a'b, \pi^{i}a'\} $$ 
for $b \in {\Cal O}_{K}$. 
So for a lifting $t_{1},...,t_{d-1}$ of a $p$-base of $F$  
we have 
$$
\aligned
\{1+\pi^{i}a', 1+\pi^{e'-i}b, t_{1},...,t_{d-1} \} & =  
\{1+ \pi^{e'}a'b, \pi^{i}, t_{1},...,t_{d-1}\} \\
& =  i\{1+ \pi^{e'}a'b, \pi, t_{1},...,t_{d-1}\} 
\endaligned
$$
in $k_{d+1}(K)$ 
(here we used $\{1+\pi^{e'}x,u_{1},...,u_{d}\}=0$ for any units $u_{i}$ 
in $k_{d+1}(K)$ 
which follows from $\Omega_{F}^{d}=0$ and the calculation of the 
subquotients $\gr_{m}k_{d+1}(K)$ in subsection 4.2). 
So we can take $b \in {\Cal O}_{K}$ such that the above symbol is non-zero 
in $k_{d+1}(K)$. 
This completes the proof in  the case where $i$ is prime to $p$. 

Next, we assume $p$ divides $i$. 
We also use the above formula, and 
calculate 
$$
\aligned
\{1+\pi^{i}a', 1+(1+b \pi) \pi^{e'-i-1}, \pi \} 
& =  \{1+ \pi^{e'-1}a'(1+b \pi), 1+ b \pi, \pi \} \\
& =  \{1+ \pi^{e'}a'b(1+b \pi), a'(1+b \pi), \pi \}. 
\endaligned
$$
Since we may think of $a'$ as a part of a lifting of a 
$p$-base of $F$, we can take some 
$x=\{1+(1+b \pi) \pi^{e'-i-1}, \pi, t_{1},...,t_{d-2} \}$ 
such that $\{a, x\} \neq 0$ in $k_{d+1}(K)$. 

\medskip

If $\ell$ is prime to $\chr(F)$, for the extension $L/K$ obtained 
by an equation $X^{\ell}=a$, we can find $x$ such that 
$\{a,x\} \neq 0$ in $K_{d+1}(K)/\ell$ in the same way as above, 
using 
$K_{d+1}(K)/\ell=\gr_{0}K_{d+1}(K)/\ell=K_{d}(F)/\ell$. 
In the case where $\chr(K)=p>0$ we can use Artin--Schreier theory 
instead of Kummer theory, and therefore we can argue in a similar
way to the previous  method.  
This completes the proof of the isomorphism theorem. 

Thus, the isomorphism theorem can be proved by computing symbols, once 
we know Bloch--Kato's theorem. 
See also a proof in \cite{K1}.

\Bib References

\rf{I} L. Illusie, Complexe de de Rham--Witt et
cohomologie cristalline,
Ann. Sci. \'Ecole Norm. Sup.(4), 12(1979), 501--661.

\rf{K1} K. Kato,  A generalization of local class field theory by 
using $K$-groups II, J. Fac. Sci. Univ. Tokyo 27 (1980), 603--683.

\rf{K2} K. Kato,  Galois cohomology of complete discrete valuation 
fields, In Algebraic $K$-theory, Lect. Notes in Math. 967, Springer
Berlin 1982, 
215--238.
  
\rf{S} J.-P. Serre,  {Corps Locaux} (third edition), 
Hermann, Paris 1968.

\endBib
 
\Coordinates

Department of Mathematics \ 
Tokyo Metropolitan University 

Minami-Osawa 1-1, Hachioji, Tokyo 192-03, Japan

E-mail: m-kuri\@comp.metro-u.ac.jp  
\endCoordinates

\vfill
\pagebreak
\end